\documentclass[12pt]{article}
\usepackage[cp1251]{inputenc}       
\usepackage{amsfonts}
\usepackage[active]{srcltx}
\usepackage{amsmath}
\usepackage{amssymb}
\usepackage{cite}
\usepackage{graphicx}
\usepackage{color}
\newenvironment{sciabstract}{%
\begin{quote} \footnotesize} 
{\end{quote}}
\newtheorem{theorem}{Theorem}

\newtheorem{proposition}{Proposition}
\newtheorem{remark}{Remark}
\usepackage{graphicx}
\title{\bf What Inverse Problems Can Be Solved Using Neural Networks?
}
\author
{A.S. Leonov$^1$, D.V. Lukyanenko$^2$\\
\\
\normalsize{$^{1}$National Research Nuclear University MEPhI,}\\
\normalsize{Moscow, Russian Federation}\\
\normalsize{$^{2}$Lomonosov Moscow State University, Faculty of Physics,}\\
\normalsize{Moscow, Russian Federation}\\
\normalsize{E-mail: asleonov@mephi.ru, lukyanenko@physics.msu.ru}
}
\date{}

\begin{document}

\maketitle

\begin{sciabstract}
The relationship between classical methods for solving inverse problems and neural network methods, which have been widely used in recent years, is examined. An answer is provided to the question of what inverse problems can be solved using neural network methods. The answer is that standard neural network algorithms can solve only conditionally well-posed problems (also referred to as well-posed according to Tikhonov). The reason for this is that the design of standard neural network algorithms violates the well-known Bakushinsky veto from the theory of ill-posed problems. It requires that any method theoretically guaranteeing a stable solution to an ill-posed problem depend explicitly on the level of data errors. The article highlights the distinction between the concept of regularization in the theory of ill-posed problems and “regularization” as understood in standard neural network training, where Bakushinsky's veto is violated. Due to this violation, any changes to the neural network architecture, weights regularization during training, or training dataset augmentation cannot alter the answer specified above. It is theoretically demonstrated that an ill-posed inverse problem can be solved only by using a parametric set of neural networks that explicitly depend on the known error levels (as parameters) of perturbed data for the inverse problem. Accordingly, as an alternative to standard neural network algorithms, \emph{neural network based regularization algorithms} (NNBRAs) are proposed to solve ill-posed inverse problems. One approach to their construction involves using \emph{a regularized dataset} to train a standard neural network (possibly without weights regularization). Such a dataset can be generated from perturbed inverse-problem data using a specific base regularization algorithm that explicitly depends on the error level of the data. Test solving the simplest finite-dimensional ill-posed inverse problem demonstrate that the accuracy of approximate solutions obtained using the NNBRA is comparable to that of the base RA. At the same time, errors of conventional neural network methods may be unacceptable.

\textbf{2010 Mathematics subject classifications}: 47A52, 47J06, 65F22, 97R40

\textbf{Keywords}: inverse ill-posed problems, regularization, neural network algorithms.
\end{sciabstract}


\section{Introduction}\label{s1}
Over the past decade, there has been unprecedented interest in using neural networks to solve various problems in science, technology, economics, and other fields. In many cases, the basis for this is the speed at which a sufficiently accurate solution to the problem can be obtained using a deep-learned neural network. This area of research is also developing rapidly in the field of solving various inverse problems.

In our article, we will focus solely on so-called \emph{regression inverse problems}, leaving aside tasks such as recognition, classification, and the like. Numerous studies of this kind exist regarding the application of neural networks in medicine, geophysics, chemical engineering, non-destructive testing, robotics, etc. (see, for example, \cite{M1,M2,G1,G2,G3,G4,Chem,NonDestr1,NonDestr2,Robots}). However, these same problems have been (or are currently being) solved using traditional methods of computational mathematics as well. Thus, there exist, in a sense, two realities: classical methods for solving inverse problems and new neural network-based methods. In the neural network literature, classical methods for solving inverse problems are often referred to as \emph{the model-driven approach}, while neural network-based methods are termed \emph{the data-driven approach}. This emphasizes that the first approach to solving the inverse problem employs a known model of the phenomenon under study (the forward problem operator), whereas the second relies solely on data to reconstruct both the governing law (the inverse problem operator), which often cannot be explicitly formalized, and the solution to the inverse problem itself.

There is an understanding within the scientific community that these methods are interconnected. Studies have been conducted that establish these connections to a certain extent (see, for example, \cite{T1,T2,T3,T4,T5}). However, in our view, these studies do not fully take into account the specific nature of inverse problems and, in particular, their frequently encountered ill-posedness. For a significant number of such problems, this manifests as the instability of solutions with respect to random perturbations in the data. It is precisely these ill-posed problems that will be considered below. In the study of the neural network approach to solving inverse problems, there is still no clear answer to the question of whether a neural network of the standard type (i.e., the kind frequently encountered in publications) can solve ill-posed (unstable) problems, or if it is suited only for solving well-posed inverse problems. In other words, does the neural network theoretically constitute a Tikhonov regularization algorithm (RA, see \cite{TA,IVT,TLY}) or not? The question of how to incorporate a priori information about the solution into neural network algorithms is also important. Furthermore, the issue of estimating the accuracy of neural network solutions for ill-posed inverse problems requires analysis.

We would like to compare classical and neural-network-based methods for solving ill-posed inverse problems in greater detail, highlighting the differences between them. This make it possible to more clearly define the scope of applicability of the latter, that is, the class of inverse problems for which neural network algorithms theoretically yield adequate solutions. This is the subject of our article.

The structure of the work is as follows. Section 2 presents the mathematical formulation of inverse problems under consideration and outlines the well-known classification of inverse problems: well-posed, ill-posed, and conditionally well-posed problems. It also notes the necessity of using the error level of approximate data when constructing regularization algorithms for solving ill-posed problems (Bakushinsky’s veto). Section 3 presents a finite-dimensional approximation scheme for the inverse problem under consideration, which is necessary for its numerical solution. It is noted that such an approximation inevitably introduces an error into the inverse problem operator during practical solution. The latter can significantly affect the accuracy of the result. Section 4 presents a mathematical model of a neural network that corresponds to the majority of practical neural network implementations for solving regression problems. The formal procedure for training a neural network and the concept of \emph{regularization (weights regularization)} adopted in neural network theory, are described. A distinction is indicated between this concept and the concept of a regularizing algorithm. Key theoretical results concerning the scope of applicability and the fundamental inapplicability of standard neural network algorithms to ill-posed (unstable) problems are proven in Section 5. Section 6 proposes a constructive approach to developing \emph{neural-network-based regularization algorithms} (NNBRAs) which, unlike conventional neural network methods, enable the solution of ill-posed (unstable) inverse problems. Section 7 presents numerical illustrations that clearly demonstrate the advantage of the proposed approach even when solving the simplest degenerate system of linear algebraic equations.

\section{Inverse problem}
By an inverse problem, we understand the solution of the operator equation
\begin{equation}\label{F}
F(x)=y,~x\in X,~y\in Y.
\end{equation}
Here $X$ and $Y$ are Banach spaces, and $F(x)$ is a known continuous operator acting from $X$ to $Y$. In the following, we assume that for a given $y \in Y$, this problem has a unique solution $\bar{x} = \bar{x}(y)$. The following classification of inverse problems is well known: well-posed, ill-posed, and conditionally well-posed. A well-posed inverse problem (in the sense of Hadamard) must, given the assumptions regarding the existence and uniqueness of the solution, possess the property of stability of the solution in $X$ with respect to perturbations of the data $y$ in the space $Y$. This means that for any exact data $y \in Y$ and any approximations $y_\delta \in Y$ given with accuracy $\delta$ (i.e., $\|y - y_\delta\| \leq \delta$), the corresponding approximate solutions $x(y_\delta)$ converge in $X$ to $\bar{x}(y)$ as $\delta \to 0$. The inverse problem may be ill-posed; that is, given the assumptions made, the solution $\bar{x}(y)$ may prove unstable with respect to perturbations of the data $y$ in $Y$. This means that for certain exact data $y \in Y$, there exist approximate data $y_\delta \in Y$ such that the approximate solutions $x(y_\delta)$ do not converge to $\bar{x}(y)$ as $\delta \to 0$. Finally, the inverse problem may be conditionally well-posed (well-posed in the sense of Tikhonov) on a well-posedness set $K \subset X$ (see \cite{Lavr,BK}, etc.). This means that its solutions belonging to the set $K$ are stable with respect to perturbations of the data if these perturbed data belong to the set $F(K) \subset Y$. Typically, the sets $K$ are characterized as compact sets in $X$ (see \cite{IVT,BK}, etc.).

The inverse problem \eqref{F} often needs to be solved on a given set $D \subset X$ representing constraints on the desired solution. These constraints are related to the available \emph{a priori information} about the solution.
To solve such an inverse problem, taking into account potential ill-posedness, one should employ Tikhonov-type regularization algorithms (RAs). This concept is defined as follows. Let the algorithm for solving the inverse problem be defined by an operator $R(y_{\delta },\delta ):~Y\times (0,\delta _{0})\rightarrow X$, used to obtain an approximate solution $x_{\delta }=R(y_{\delta },\delta )$.
For each $x\in D$, we introduce the quantity \[\Delta (x,\delta )=\underset{y_{\delta }\in Y}{\sup }\left\{ \left\Vert x-R(y_{\delta },\delta )\right\Vert :~\left\Vert F(x)-y_{\delta }\right\Vert \leq \delta ~\right\}.\] It represents the worst-case accuracy of the approximate solution to equation \eqref{F}, obtained using algorithm $R$, for a fixed solution $x$ and all possible admissible approximate data. Algorithm $R$ is called a regularizing algorithm for equation (1) on the class of solutions $D$ if $\lim_{\delta \to 0} \Delta(x, \delta) = 0$ for all $x \in D$. \emph{Solving an ill-posed problem} means specifying a particular regularization algorithm for it and applying that algorithm in practice. Clearly, such an algorithm is also suitable for solving well-posed and conditionally well-posed problems. It should be noted that the explicit dependence of the regularization algorithm on the data errors level $\delta$ is fundamental. As proven in \cite{BakGon}, it is impossible to solve an ill-posed inverse problem without knowledge and explicit use of the value of $\delta$. Methods of the form $x_{\delta }=R(y_{\delta })$ that do not explicitly depend on $\delta$ are suitable only for solving well-posed or conditionally well-posed problems \cite{BK}. This constitutes the so-called Bakushinsky veto.

The above holds true also when the spaces $X$ and $Y$ are finite-dimensional. Ill-posed problems of this type also quite frequently as well. Characteristic example is computing the pseudoinverse of a singular or rank-deficient matrix, as well as finding the normal solution to a singular system of linear algebraic equations. These problems will be discussed in greater detail below.

\section{Finite-dimensional approximation of the inverse problem}
Solving inverse problems in infinite-dimensional spaces typically requires their finite-dimensional approximation. For the sake of definiteness, we consider one of the existing forms of such approximation.

Let a sequence of finite-dimensional spaces $X_n$, $\mathrm{dim} X_n=n$, and a sequence of finite-dimensional sets $D_n$ ($D_n \subset X_n$) be defined, such that: \vspace{-2mm}
\begin{equation*}
X_{1}\subset X_{2}\subset ...\subset X_{n}\subset ...\subset X,\,\
\bigcup_{n=1}^{\infty }X_{n}=X,\,\, {D}_{1}\subset
{D}_{2}\subset ...\subset {D}_{n}\subset
...\subset {D},\,\ \bigcup_{n=1}^{\infty
}{D}_{n}={D}. \vspace{-2mm}
\end{equation*}
For each $n$, we assume given an operator $F_{n}: {D}_{n} \rightarrow Y$ continuous on $D$ and an element $y_{n} \in Y$ that approximate the data $(F, y)$ of problem (\ref{F})
in the following sense:
$\left\| F_{n}(x)-F(x)\right\| \leq \psi(h_{n},\left\| x\right\|) \quad \forall x\in {D}$ and $\left\| y_{n}-{y}\right\| \leq \delta _{n}$,
where $h_{n}, \delta _{n} \rightarrow 0$ as $n \rightarrow \infty$.
Typical examples of approximation measures $\psi(h, \|x\|)$ for the operator  can be found, for instance, in \cite{TLY}. Thus, the approximate data for the inverse problem are given by the quantities $S_n = (F_n, y_n, h_n, \delta_n)$. In other words, these represent approximate realizations of the operator and its right-hand side in equation \eqref{F}, along with error estimates for this data. In this case, the inverse problem itself consists of finding, for each given $n$ and given data $S_n=(F_{n},y_{n},h_{n},\delta_{n})$, an approximation $x_n(S_n)\in X_n$ to the element $\bar{x}(y)$. The method for finding such an approximation is a regularization algorithm if, for any data $S_n$ of the specified type, convergence $x_n(S_n)\overset{X}{\rightarrow }\bar{x}(y)$ is ensured as $n\to\infty$. Examples of such RAs can be found, for instance, in \cite{TLY,EHN}, among others.

\section{Neural network model}\label{sec4}
We will solve equations (1) on the set $D$ using neural networks. Their application is typically based on the use of finite-dimensional discrete quantities. However, we first consider neural network models in the general case, without assuming the finite-dimensionality of the spaces $X$ and $Y$ or the set $D \subset X$. Consequently, the finite-dimensional special case is encompassed by this treatment; its specific features will be discussed separately. Let us specify what we will refer to as a neural network (NN) in this article. We will not address concepts such as NN architecture at this stage. Instead, we will examine potential structures and types of training sets, loss functions, and training procedures separately, as the need arises. We consider an untrained neural network to be a fixed parametric family of continuous operators $NN_{w}:D\rightarrow X$ depending on parameters $w=\{w_{1},...,w_{M}\}$. To train the network, we use a dataset $\mathcal{D}=\{R_{X},D_{Y}\}$, where $R_{X}\subset D$ is a set of solutions to the inverse problem and $D_{Y}\subset F(D)$ is the set of corresponding data. Thus, $D_{Y}=F(R_X)$ and ${\cal D} = \{ (x,y):\,\,x \in {R_X},\,y = F(x) \in {D_Y}\} $. Let us specifically single out subsets $R_{X}^{\prime }\subset R_{X}$ and $D_{Y}^{\prime }\subset D_{Y}$ such that $F(R'_{X})=D'_{Y}$. Accordingly, the set $\mathcal{D}'=\{R_{X}^{\prime },D_{Y}^{\prime }\}$ is interpreted as the training subset, while the set $\mathcal{D}\backslash \mathcal{D}'$ constitutes the network validation set. Training consists of finding the parameters $\bar{w}$ as the solution to the variational problem
\begin{equation}\label{NN}
\bar{w}=\arg \min \left\{ \Phi \lbrack R_{X}^{\prime },NN_{w}(D_{Y}^{\prime })]:~w\in W\right\},
\end{equation}
where $\Phi (\cdot,\cdot) $ is a given continuous function representing the aggregate residual between the solutions from the set $R_{X}^{\prime }$ and the neural network responses, for a fixed $w$, to the corresponding data $D_{Y}^{\prime }$. The set $W$ is defined by possible constraints on the weight coefficients. A characteristic example of the function $\Phi$ is the maximum deviation of the solutions in the dataset from the neural network's solutions for a given $w$:
\begin{equation}\label{Loss}
\Phi [R_X^\prime ,N{N_w}(D_Y^\prime )] = \max \left\{ {\rho [x',N{N_w}(y')]:{\mkern 1mu} \,\forall {\mkern 1mu} x' \in R_X^\prime ,y'=F(x') \in D_Y^\prime } \right\} = f(w),
\end{equation}
where $\rho [x',N{N_w}(y')] = \left\| {x' - N{N_w}(y')} \right\|$ is \emph{the loss function}. A metric defined on $X$ can also be used as the loss function. Special optimization methods have been developed to solve problem \eqref{NN} for finite-dimensional spaces $X$ and $Y$, but we will not discuss them here.

We will briefly denote the trained neural network $NN_{\bar{w}}$ with loss function $\rho$ as $NN$, and the approximate solution obtained using the neural network for data $y \in D$ as $x_{NN} = NN(y)$. Important characteristics of this neural network are the training accuracy
\[\varepsilon _{0}=\min \left\{ \Phi \lbrack R_{X}^{\prime },NN_{w}(D_{Y}^{\prime })]:~w\in W\right\} =\Phi \lbrack R_{X}^{\prime },NN(D_{Y}^{\prime })]=f(\bar w)\]
and the validation accuracy:
\[\varepsilon _{val}=\max \left\{ \left\Vert x-NN(y)\right\Vert :~\forall (x,y)\in \mathcal{D}\backslash \mathcal{D}'\right\}.\]
The accuracy of the neural network's performance over the set of all admissible data is of particular importance:
\begin{equation}\label{Toch}
\varepsilon =\max \left\{ \left\Vert F^{-1}(y)-NN(y)\right\Vert :~\forall y\in F(D)\right\}.
\end{equation}
The operational algorithm of many neural networks corresponds to the theoretical scheme described.

Above we made an assumption about the continuity of the operator $NN_w$. This feature is typical for many neural networks. As an example, consider popular finite-dimensional networks such as a multilayer perceptron with $N$ layers. In them, the connection between the input $\{x_{j}^{(l-1)}\}$ and the output $\{x_{j}^{(l)}\}$ of layer number $l$ is determined by a formula like
\begin{equation}\label{NN1}
x_j^{(l)} = f\left( \sum_{i=1}^{n_{l-1}} w_{j,i}^{(l)} x_i^{(l-1)} + b_j^{(l)}\right),\, j=1,...,n_{l}.
\end{equation}
Here, $f(\cdot)$ is an activation function (typically continuous), $w_{j,i}^{(l)}$ represents the weight matrix of the network's $l$-th layer, $b_{j}^{(l)}$ are the bias coefficients, and $n_{l-1}$ is the fixed number of neurons in layer $(l-1)$ (see, e.g., \cite{Pers}). By setting $\{x_{j}^{(0)}=y\}$ and $x_{NN}=\{x_{j}^{(N)}\}$, and taking into account the architecture of such a neural network, one can observe the continuity of the operator $NN_w$ as a superposition of expressions of the form \eqref{NN1}.

The variational problem \eqref{NN} often turns out to be ill-posed. Their ill-posedness may manifest as non-uniqueness of the solution, ill-conditioning of the extremum point, or its instability with respect to perturbations of the data $\mathcal{D}'$. Theoretically, such a variational problem can be solved using Tikhonov regularization (see, e.g., \cite{TLY}).
However, special iterative methods are typically employed when training neural networks in general, and specifically those when problem \eqref{NN} is ill-posed. Their implementation can lead to overtraining, a phenomenon in which the objective function may begin to increase from a certain iteration onward during the minimization process. One way to mitigate this effect is through weight regularization, which is interpreted in neural network theory as follows. The objective function $\Phi \lbrack D_{X}^{\prime },NN_{w}(D_{Y}^{\prime })]$ is penalized by a strictly convex functional $\Omega \lbrack w]$ of the weights, that is, problem \eqref{NN} is replaced by the following:
\begin{equation}\label{Loss1}
{w_\lambda } = \arg \min {\rm{ }}\left\{ {\Phi [D_X^\prime ,N{N_w}(D_Y^\prime )] + \lambda \Omega [w]:w \in W} \right\}.
\end{equation}
Using the found $w_\lambda$, the neural network computes an approximate solution $x_{\delta\lambda} = NN_{w_\lambda}(y_\delta)$ for approximations $y_\delta$ to the exact data $y$.
Here, $\lambda > 0$ is an analogue of the regularization parameter used when solving the variational problem via Tikhonov's method \cite{TA,TLY}. It should be noted that in the neural network approach, $\lambda$ is typically chosen heuristically and does not explicitly depend on the data error level $\delta$ of the inverse problem being solved. Thus, this procedure disregards Bakushinsky’s veto, and such a method cannot be considered a regularization algorithm for solving the inverse problem. This means that it does not guarantee the convergence of the approximations $x_{\delta\lambda}$ to the exact solution $\bar{x}(y)$ as $\delta \to 0$.

\section{Key assertions}
\begin{proposition}\label{p1}
Let problem (1) be conditionally well-posed on a compact set $K, K \subset X$. Then a neural network can solve the problem with an accuracy no better than $\bar{\varepsilon}=\max \left\{\varepsilon_0,\varepsilon _{val}\right\}$.
\end{proposition}
Proof. Conditional well-posedness, that is the stability of the inverse problem with respect to data perturbations that do not take the data outside the set $F(K)$, implies the continuity of the operator $F^{-1}(y)$ on $F(K)$ \cite{Lavr,BK}. Therefore, due to the continuity of the operator $NN(y)$ on $F(K)$, the functional $\left\Vert F^{-1}(y)-NN(y)\right\Vert$ is continuous on the set $F(K)$. Then the accuracy \eqref{Toch} of the neural network on the set $K$ is finite, and the following estimate holds:
\begin{multline*}
\varepsilon =\max \left\{ \left\Vert F^{-1}(y)-NN(y)\right\Vert :~\forall y\in F(K)\right\} \geq  \\\ge \max \left\{ \left\Vert x-NN(y)\right\Vert :~\forall x\in R_{X},y\in D_{Y}\right\}\ge \max \left\{\varepsilon_0,\varepsilon _{val}\right\}=\bar{\varepsilon},
\end{multline*}
which is obtained taking into account the inclusions ${R_X} \subset K$ and ${D_Y} \subset F(K)$. $\square$

Since the value of $\varepsilon$ from \eqref{Toch} is typically non-zero, the inverse problem on a general set $D$ should be solved using sequences of neural networks with a decreasing value of $\varepsilon$. This can be achieved, for example, by expanding datasets, increasing the dimensionality of the finite-dimensional spaces $X_n$, etc. The following statement concerns such sequences of neural networks and serves, to a certain extent, as a converse to Proposition 1.
\begin{theorem}\label{p2}
Let the set $D$, as a subset of $X$, be dense in itself.
Furthermore, let there be a sequence of neural networks $\{NN^{(n)}\}_{n=1}^{\infty}:\, Y\to D$ solving the inverse problem on the set $D$ with individual accuracy $\varepsilon_n > 0$ (see \eqref{Toch}), where $\varepsilon_n\to 0$ as $ n\to \infty$. Then the inverse problem being solved is conditionally well-posed on $D$.
\end{theorem}
Proof. Let $y \in F(D)$ be an arbitrary element. Since the set $D$ is dense-in-itself and the operator $F$ is continuous, the set $F(D)$ is also dense-in-itself. Therefore, there exists a sequence $\{y_{m}\}_{m=1}^{\infty }\subset F(D)$ such that $y_{m}\rightarrow y$ as $m\rightarrow \infty $. Let us take an arbitrary sequence $\{y_{m}\}_{m=1}^{\infty }$ of this kind. For a fixed $n$, the following estimate holds:
\begin{multline*}
\left\| {{F^{ - 1}}(y) - {F^{ - 1}}({y_m})} \right\| \le \left\| {{F^{ - 1}}(y) - N{N_n}(y)} \right\| + \left\| {N{N_n}(y) - N{N_n}({y_m})} \right\| \\ + \left\| {N{N_n}({y_m}) - {F^{ - 1}}({y_m})} \right\|\le
\varepsilon_n +\left\| {N{N_n}(y) - N{N_n}({y_m})} \right\|+\varepsilon_n,
\end{multline*}
which follows from \eqref{Toch}. From this, passing to the limit as $m \to \infty$, we obtain:
\[\mathop {\lim }\limits_{m \to \infty } \left\| {{F^{ - 1}}(y) - {F^{ - 1}}({y_m})} \right\| \le 2{\varepsilon _n} + \mathop {\lim }\limits_{m \to \infty } \left\| {N{N_n}(y) - N{N_n}({y_m})} \right\| = 2{\varepsilon _n}.\]
Passing to the limit as $n \to \infty$ in this estimate, we arrive at the equality \[\mathop {\lim }\limits_{m \to \infty } \left\| {{F^{ - 1}}(y) - {F^{ - 1}}({y_m})} \right\| = 0.\] This means that the operator $F^{-1}$ is relatively continuous on the set $D$, and thus the inverse problem under consideration is conditionally well-posed on $D$. $\square$

Thus, an inverse problem can be solved with any arbitrarily high degree of precision using neural networks only if it is conditionally well-posed on the set $D$.

\begin{theorem}\label{p3}
Let problem (1) be ill-posed (unstable with respect to data perturbations in $Y$), and for some element $\hat{y}$ from the dataset there is a sequence $\{ {y_n}\}$ such that ${y_n} \in D\backslash {D_Y},\,{y_n} \to \hat y $, and for it the inequality holds
\begin{equation}\label{eq3}
\left\| {{F^{ - 1}}({y_n}) - {F^{ - 1}}(\hat y)} \right\| \ge E = const >\max \left\{\varepsilon_0,\varepsilon _{val}\right\}=\bar{\varepsilon}.
\end{equation}
Then the accuracy of solving the problem using a neural network cannot be better than $E-\bar{\varepsilon}$.
\end{theorem}
Proof. For the element $\hat{y} \in D_Y$ and the sequence $\{y_n\}$ specified in the theorem's hypotheses, the following inequalities hold:
\[\left\Vert F^{-1}(y_{n })-NN(y_{n})\right\Vert \ge \left\Vert F^{-1}(y_{n })-F^{-1}(\hat{y})\right\Vert -\left\Vert F^{-1}(\hat{y})-NN(\hat{y})\right\Vert -\left\Vert NN(\hat{y})-NN(y_{n })\right\Vert.\] Hence, by condition \eqref{eq3}, by inequality of training and validation
\[\left\Vert F^{-1}(\hat{y})-NN(\hat y)\right\Vert \leq \max \left\{\varepsilon_0,\varepsilon _{val}\right\}=\bar{\varepsilon}\]
and taking into account the continuity of the operator $NN$, i.e. equality
\[\underset{n \rightarrow \infty}{\lim }\left\Vert NN(\hat{y})-NN(y_{n })\right\Vert =0,\]
we get
\[\underset{n\rightarrow \infty }{\underline{\lim }}\left\Vert F^{-1}(y_{n})-NN(y_{n})\right\Vert \geq E-\bar{\varepsilon}.\]
This means that when processing data $\{y_{n}\}$ that are sufficiently close to the sample points $\hat{y} \in D_{Y}$ but do not belong to the dataset, the neural network yields approximate solutions $NN(y_{n})$ with an accuracy no better than $E-\bar{\varepsilon}$. Thus, its overall performance accuracy \eqref{Toch} is also no better than $E-\bar{\varepsilon}$. $\square$

Inequality \eqref{eq3} expresses the instability of computing the element ${F^{- 1}}(\hat y)$, reflecting the ill-posed nature of the inverse problem. Theoretically, the quantity $E$ can be arbitrarily large.
This can lead to arbitrarily large errors in the neural network's solution for elements not included in the dataset. Thus, neural networks can solve only well-posed or conditionally well-posed inverse problems.
This also applies to the case of finite-dimensional spaces $X$ and $Y$. A corresponding numerical illustration is provided in Section \ref{Ill}.

In the case of finite-dimensional spaces, the boundedness of a closed set $D$ implies its compactness, which guarantees conditional well-posedness. Thus, finite-dimensional neural network algorithms are capable of solving the inverse problem on such sets $D$, provided that the data do not lie outside the set $F(D)$.

\section{How to solve ill-posed problems using neural networks}\label{s7}
To enable neural networks to solve ill-posed problems, it is necessary to take Bakushinsky's veto into account during their design. This means that neural networks must be applied to inverse problems involving approximate data for which an accuracy estimate $\delta$ is known. Furthermore, the values of $\delta$ must be taken into account during the actual process of creating the neural network. This can be achieved, for example, by constructing the dataset not from pairs of the form $(x, y=F(x))$, but from pairs $(x=R_\delta(y), y)$, where $R_\delta(\cdot) = R(\cdot, \delta)$ is a certain RA. We will refer to a dataset of this type as regularized, and the method of solving the inverse problem using such neural networks as \emph{neural-network-based regularization} (NNBR). The regularizing properties of such neural networks are realized as follows.

Let problem \eqref{F} have a solution $\bar{x} \in D$ corresponding to the data $\bar{y} = F(\bar{x})$. We solve it for approximate data $y_{\delta} \in Y$, with $\left\Vert \bar{y} - y_{\delta} \right\Vert \leq \delta$, using a regularized neural network $NN_{\delta}$ trained on the dataset $\{(x = R_\delta(y), y)\}\subset X\times Y$. Assume that for the data $y_{\delta}$, there exists a dataset element $\hat{y}_{\delta}$ such that $\left\Vert \hat{y}_{\delta} - y_{\delta} \right\Vert \leq \delta$. Then, the neural network provides the appropriate solution to the problem, $\hat{x}_{\delta }=NN_{\delta }(\hat{y}_{\delta })$, with an accuracy $\Delta (\delta )$ such that
\begin{equation}\label{accur}
\Delta (\delta ) = \left\| {\bar x - {{\hat x}_\delta }} \right\| = \left\| {\bar x - {R_\delta }({{\hat y}_\delta })} \right\| \le \left\| {\bar x - {R_\delta }({y_\delta })} \right\| + \left\| {{R_\delta }({y_\delta }) - {R_\delta }({{\hat y}_\delta })} \right\|.
\end{equation}

We consider now a sequence of such neural networks, $NN_{\delta _{n}}$, with $\delta _{n}\rightarrow 0$.
By definition of RA $R_{\delta _{n}}(y_{\delta _{n}})\rightarrow \bar{x},R_{\delta _{n}}(\hat{y}_{\delta _{n}})\rightarrow \bar{x}$, because $y_{\delta _{n}},\hat{y}_{\delta _{n}}\rightarrow \bar{y}$ as $n\rightarrow \infty $. Therefore, from estimate \eqref{accur}, we obtain
$\Delta (\delta _{n})=\left\Vert \bar{x}-\hat{x}_{\delta _{n}}\right\Vert \rightarrow 0$, i.e. a sequence of described neural networks provides regularization of the inverse problem on an arbitrary element $\bar{x}\in D$.

The solution of linear ill-posed problems using neural network regularization algorithms is of particular interest. Let $A:X\rightarrow Y$ be a bounded linear operator, and let us consider the equation $Ax=y$ on a set $D$. In the construction of NNBR networks, we employ Tikhonov regularization: $x^{\alpha (\delta )}=R_{\delta }(y_{\delta })=(\alpha (\delta )I+A^{\ast }A)^{-1}A^{\ast }y_{\delta }$ with the choice of the regularization parameter, which will be discussed below. Then the accuracy \eqref{accur} for such neural networks $NN_{\delta }$ can be estimated as follows:
\begin{multline*}
\Delta (\delta )=\left\Vert \bar{x}-\hat{x}_{\delta }\right\Vert \leq \left\Vert \bar{x}-(\alpha (\delta )I+A^{\ast }A)^{-1}A^{\ast }y_{\delta }\right\Vert +\\
\left\Vert (\alpha (\delta )I+A^{\ast }A)^{-1}A^{\ast }y_{\delta }-(\alpha (\delta )I+A^{\ast }A)^{-1}A^{\ast }\hat{y}_{\delta }\right\Vert \\
\leq \left\Vert \bar{x}-(\alpha (\delta )I+A^{\ast }A)^{-1}A^{\ast }\bar{y}\right\Vert +\left\Vert (\alpha (\delta )I+A^{\ast }A)^{-1}A^{\ast }(\bar{y}-y_{\delta })\right\Vert\\
+\left\Vert (\alpha (\delta )I+A^{\ast }A)^{-1}A^{\ast }(y_{\delta }-\hat{y}_{\delta })\right\Vert \\
\leq \left\Vert \bar{x}-(\alpha (\delta )I+A^{\ast }A)^{-1}A^{\ast }A\bar{x}\right\Vert +\frac{\left( \left\Vert \bar{y}-y_{\delta }\right\Vert +\left\Vert \hat{y}_{\delta }-y_{\delta }\right\Vert \right) }{\sqrt{\alpha (\delta )}}\\
\leq \left\Vert \alpha (\delta )(\alpha (\delta )I+A^{\ast }A)^{-1}\bar{x}\right\Vert +\frac{2\delta }{\sqrt{\alpha (\delta )}}=
\omega _{1}(\alpha (\delta ))+\frac{2\delta }{\sqrt{\alpha (\delta )}}.
\end{multline*}
Standard techniques \cite{M,Vai,EHN} were used to obtain this estimate.
It is known that $\omega _{1}(\alpha (\delta ))\rightarrow 0$ if $\alpha (\delta )\rightarrow 0$ as $\delta \rightarrow 0$ (see \cite{M,Vai,EHN}). Therefore, if the regularization parameter is chosen such that $\alpha (\delta )\rightarrow 0$ and $\frac{\delta }{\sqrt{\alpha (\delta )}}\rightarrow 0$, then convergence $\Delta (\delta )\rightarrow 0$ as $\delta \rightarrow 0$ is guaranteed. This implies that the family of neural networks $NN_{\delta _{n}}$ constitutes a neural-network based regularization algorithm for solving the problem under consideration.

A similar result can be obtained when the regularization parameter is chosen according to the discrepancy principle \cite{LeoBook}. If a priori information regarding the source representability of the solution $\bar{x}$ via the operator $A^{\ast }A$ is available, then \emph{a priori accuracy estimates} for neural network algorithms can be derived by analogy with \cite{IVT,BK,Vai,EHN}. It should also be noted that other known regularization methods (see, e.g., \cite{IVT,BK,Vai,EHN}, etc.) can be employed as the RA instead of Tikhonov regularization.
\begin{remark}
An analogue of the concept of a regularized dataset was previously considered in \cite{T5}. However, in its construction, the authors did not take into account Bakushinsky's veto and the explicit dependence of such a dataset on the data error level $\delta$.
\end{remark}

\section{Illustrations}\label{Ill}
\subsection{Problem of computing the pseudoinverse matrix for a singular system of linear algebraic equations.}
Let us consider one of the simplest ill-posed problems. Suppose we need to construct a neural network to find the normal pseudo-solution of the system of linear algebraic equations
\[\left\{ \begin{array}{l}
{a_{11}}x + {a_{12}}y = {b_1} \\
{a_{21}}x + {a_{22}}y = {b_2} \\
\end{array} \right.\]
with a potentially singular matrix $A=[a_{i,j}]_{i,j=1}^{2,2}$ of rank at least 1. To clarify the terminology, we recall that a pseudo-solution of a system of linear algebraic equations is any element that minimizes the system's residual. By a normal pseudo-solution, we mean a pseudo-solution with minimal 2-norm.
The normal pseudo-solution can be calculated using the pseudoinverse matrix $A^+ = \lim_{\alpha \to +0} (\alpha I + A^T A)^{-1} A^T$ according to the rule $\bar{X}=(\bar{x},\bar{y})^{T}=A^{+}(b_{1},b_{2})^{T}$. Simple examples demonstrate the instability of the procedure for computing the pseudoinverse matrix and, consequently, the instability of the normal pseudo-solution with respect to perturbations of singular matrix $A$. Indeed, let
\[A = {\mathop{\rm diag}\nolimits} \left( {1,0} \right),\,{A_h} = {\mathop{\rm diag}\nolimits} \left( {1,h} \right),\,\left\| {A - {A_h}} \right\|_2 = h\to 0.\]Then \[{A^ + } = {\mathop{\rm diag}\nolimits} \left( {1,0} \right),\,A_h^ + = {\mathop{\rm diag}\nolimits} \left( {1,\frac{1}{h}} \right),\,{\left\| {{A^ + } - A_h^ + } \right\|_2} = \frac{1}{h}\to \infty.\]
Suppose we use a neural network to solve the problem for given $h$, and its dataset consists of pairs of the form $\left( A_h^ + ,A_h \right)$. Then, in this example, condition \eqref{eq3} can be written as follows:
\[\left\| {{F^{ - 1}}({y_n}) - {F^{ - 1}}(\hat y)} \right\| = {\left\| { A_{h_n}^+ -A^+} \right\|_2} = \frac{1}{h_n} = E,\]where $h_n\to 0$.
Consequently, by Theorem \ref{p3}, a neural network, regardless of its training and validation accuracy, cannot solve such a problem with adequate precision. Moreover, including increasingly accurate approximations $A_{h_n}$ of $A$ in the dataset, that is, reducing the value of $h_n$ does not alter the situation.

However, according to the results of Section \ref{s7}, the use of the regularized dataset $\left\{ {\left( {{{\left( {\alpha (h)I + A_h^ * {A_h}} \right)}^{ - 1}}A_h^ * ,{A_h}} \right)} \right\}$ guarantees obtaining a neural-network regularized approximation to the pseudoinverse matrix $\bar{A}^{+}$.

\subsection{Numerical experiments}
The statements made above can be illustrated by numerical experiments. Here we present one such experiment. We seek the normal solution $X=(x,y)^T=(1,1)^T$ to the degenerate model system of linear algebraic equations
\[\left\{ \begin{array}{l}
x + y = 2 \\
2x + 2y = 4 \\
\end{array} \right.\]
using the exact right-hand side $(b_1,b_2)^T=(2,4)^T$ and the perturbed matrix
\[{A_h} = A + H = \left( {\begin{array}{*{20}{c}}
1 & 1  \\
2 & 2  \\
\end{array}} \right) + \left( {\begin{array}{*{20}{c}}
{{h_{11}}} & {{h_{12}}}  \\
{{h_{21}}} & {{h_{22}}}  \\
\end{array}} \right),{\mkern 1mu} {\mkern 1mu} {\left\| H \right\|^2} = \sum\limits_{i = 1}^2 {\sum\limits_{j = 1}^2 {h_{ij}^2} } ,\,\,h = \frac{{\left\| H \right\|}}{{\left\| A \right\|}}.\]
The range of the operator $A + H$ may turn out to be larger than the range of the operator $A$, and therefore the problem under consideration is ill-posed. In our experiments, the right-hand side of the system of linear algebraic equations was not perturbed for ease of analyzing the results.

The problem was solved as follows. First, using the data $A_h$ and one of the algorithms under consideration, an approximation $A(h)^+$ to the pseudoinverse matrix $A^+$ was obtained. Then, this approximation was used to compute the approximate solution $X_h = A(h)^+(b_1, b_2)^T$ to the system of linear equations and its relative accuracy $\Delta_h = \frac{{{{\left\| {{X_h} - X} \right\|}_2}}}{{{{\left\| X \right\|}_2}}}$. In this context, several solution approaches are compared: a standard regularization algorithm $R_h(A_h)$ for the approximate computation of the matrix $A^+$; a NNBR algorithm trained on a regularized dataset obtained using this RA; and a neural network algorithm trained on a standard (non-regularized) dataset of the form $\left\{ {\left( {A_h^ + ,{A_h}} \right)} \right\}$. In the latter case, the matrices $A_h^+$ were calculated analytically using the formula for the pseudoinverse of a $2 \times 2$ matrix, which follows from the well-known equality $A_h^+ = \lim_{\alpha \to +0} (\alpha I + A_h^T A_h)^{-1} A_h^T$.

\begin{figure}[t]
\centering
\includegraphics[width=1\textwidth]{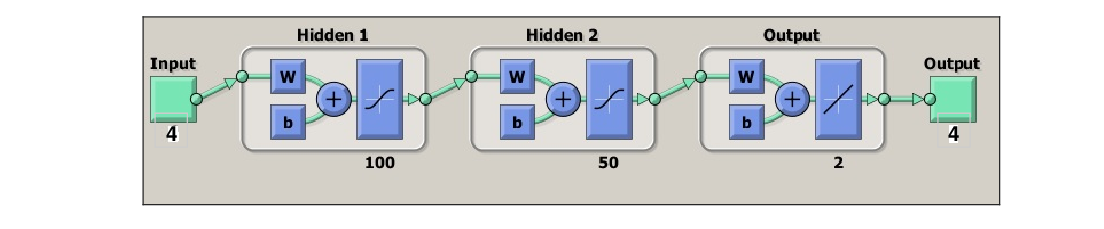}
\caption{Architecture of the neural networks used. The numbers shown indicate the dimensions of the input and output data, as well as the number of neurons in the layers.}
\label{F0}
\end{figure}

\begin{figure}[t]
\centering
\includegraphics[width=1\textwidth]{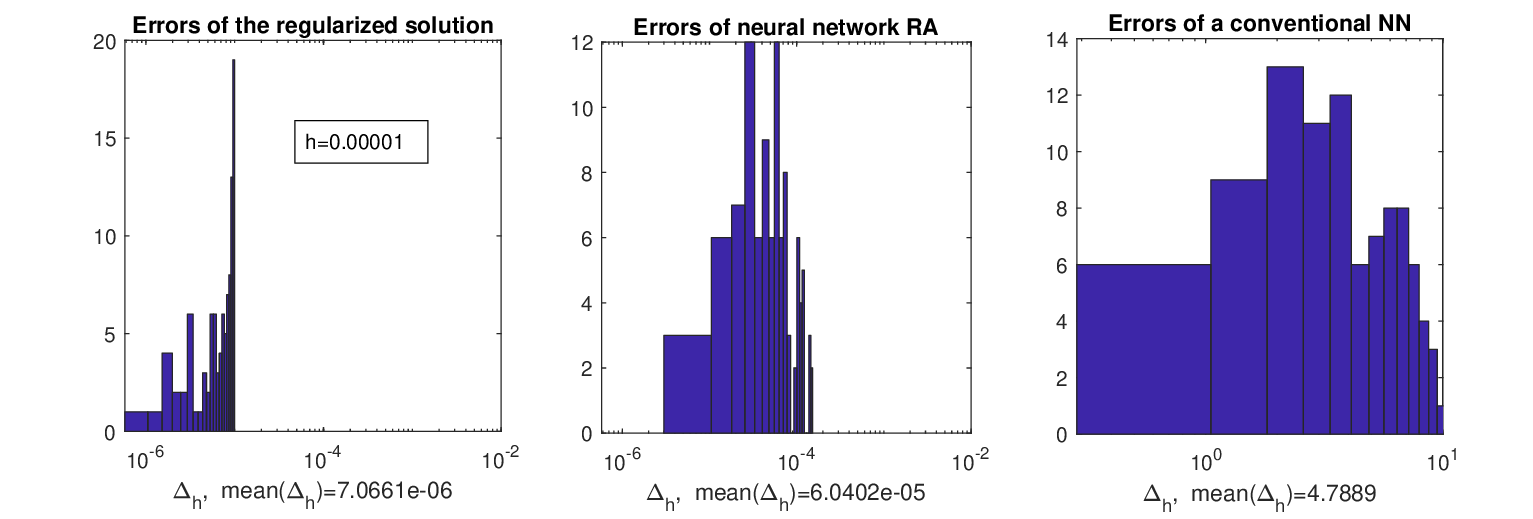}%
\caption{The distribution of errors in the three methods under consideration for $h=0.00001$.}
\label{F1}
\end{figure}

\begin{figure}[t]
\centering
\includegraphics[width=1\textwidth]{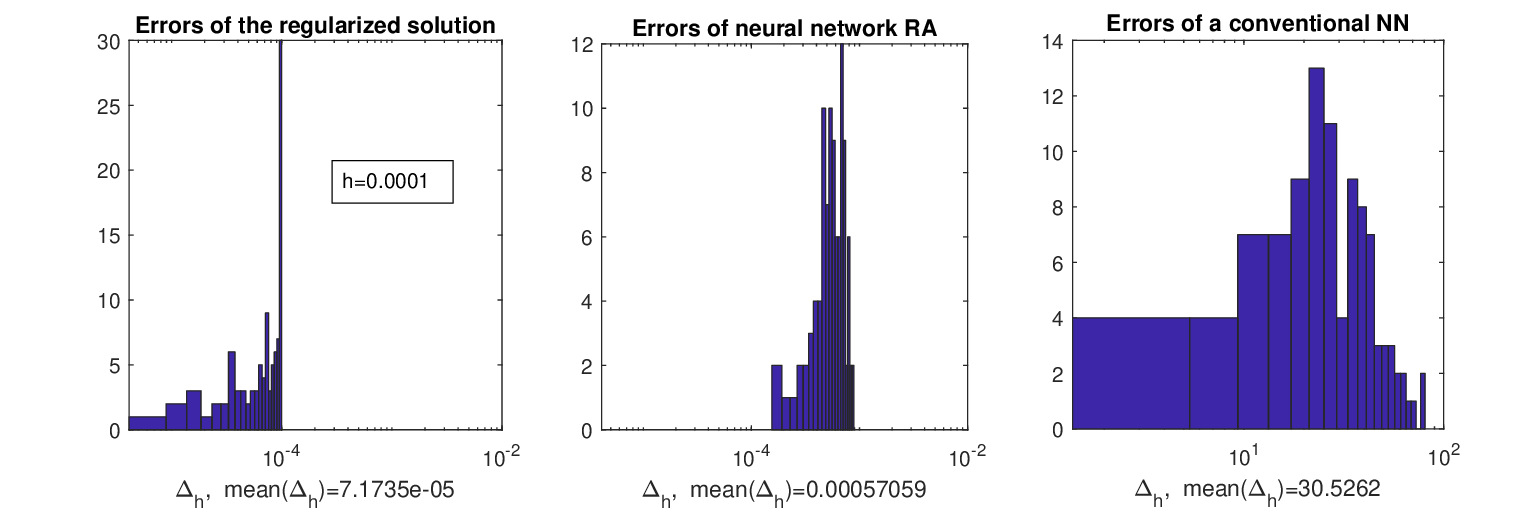}%
\caption{Distribution of errors in the three methods under consideration for $h=0.0001$.}
\label{F2}
\end{figure}

\begin{figure}[t]
\centering
\includegraphics[width=1\textwidth]{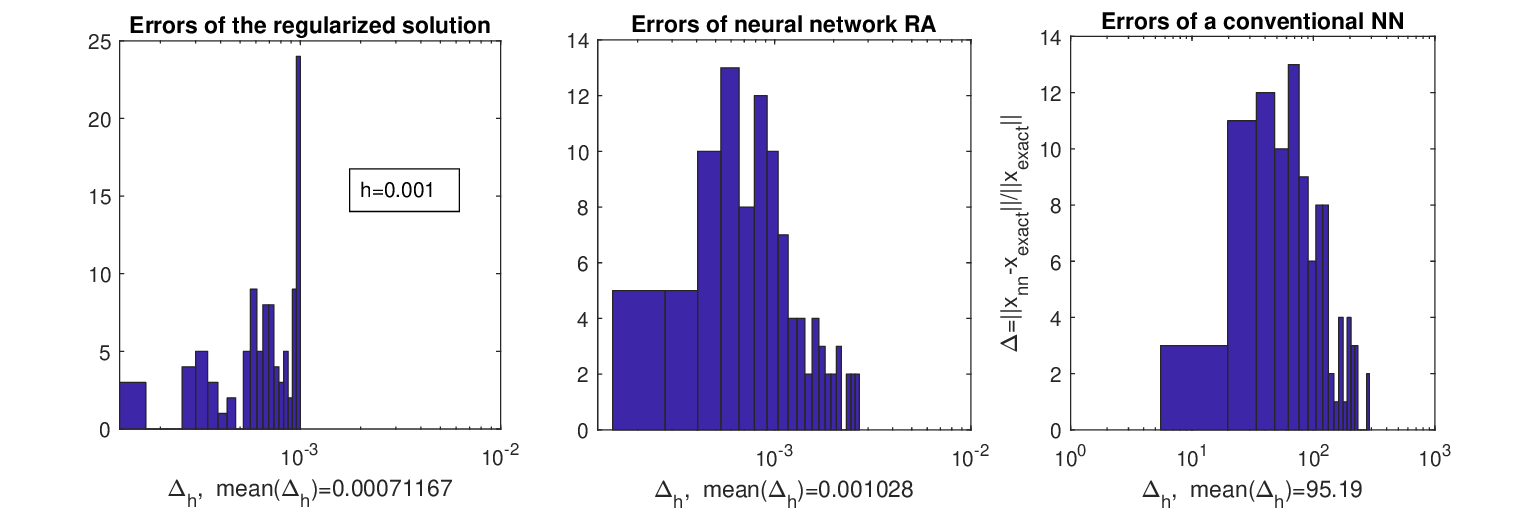}%
\caption{Distribution of errors in the three methods under consideration for $h=0.001$.}
\label{F3}
\end{figure}

\begin{figure}[t]
\centering
\includegraphics[width=1\textwidth]{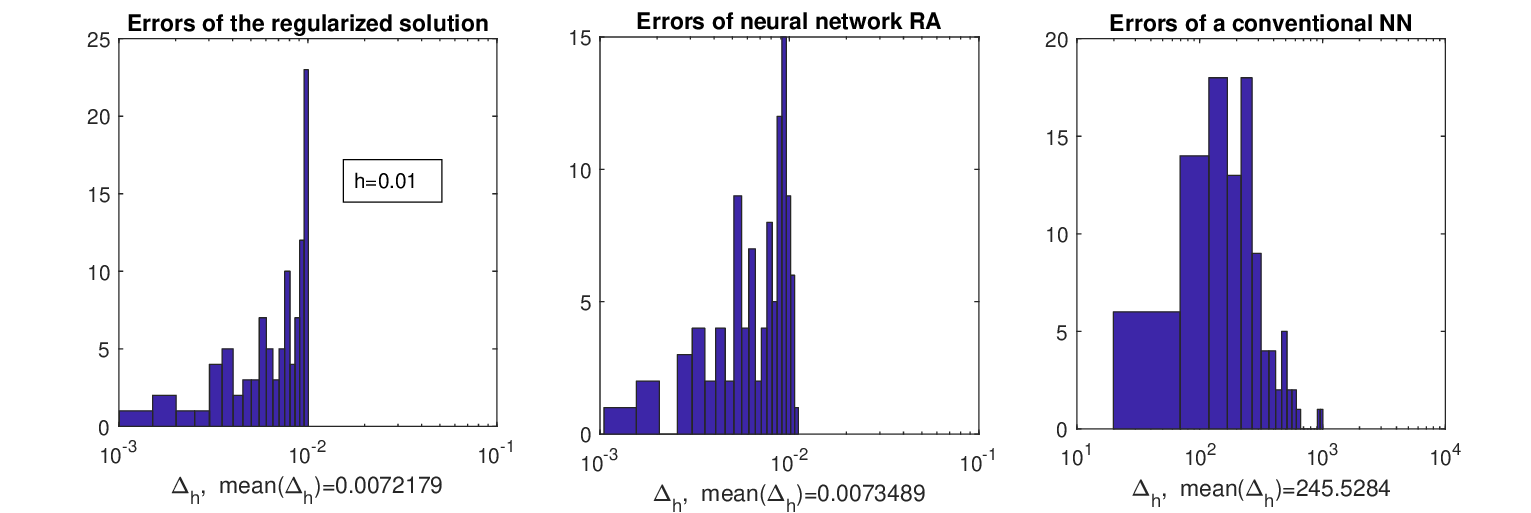}%
\caption{Distribution of errors in the three methods under consideration for $h=0.01$.}
\label{F4}
\end{figure}

The simplest variant of the TSVD method (see, e.g., \cite{EHN}) was used as the regularization algorithm $R_h(A_h)$. This method uses the singular value decomposition ${A_h} = {U_h}{\mathop{\rm diag}\nolimits} ({s_1}(h),{s_2}(h))V_h^T$, ${{s_1(h)}} \ge {{s_2(h)}} \ge {{0}}$. Regularization in the calculation of $A_h^+$ consists of replacing the singular value $s_2(h)$ with zero when the condition ${s_2}(h)<h^2$ is met, such that ${R_h}({A_h}) = {V_h}{\rm{diag}}(s_1^{ - 1}(h),0)U_h^T$.
Neural network solutions were obtained with the help of fully connected multilayer perceptrons with the architecture shown in Fig.\ref{F0}.
Training datasets were generated from matrices of the form $A + H$ for each considered perturbation level $h$ and consisted of 50,000 random samples. During testing, the methods were compared using a large number of random random matrices $H$ (50,000 realizations). Figures \ref{F1}–\ref{F4} show histograms of the accuracy distributions for the approximate solutions obtained using the specified methods at various levels of $h$.

\begin{figure}[t]
\centering
\includegraphics[width=1\textwidth]{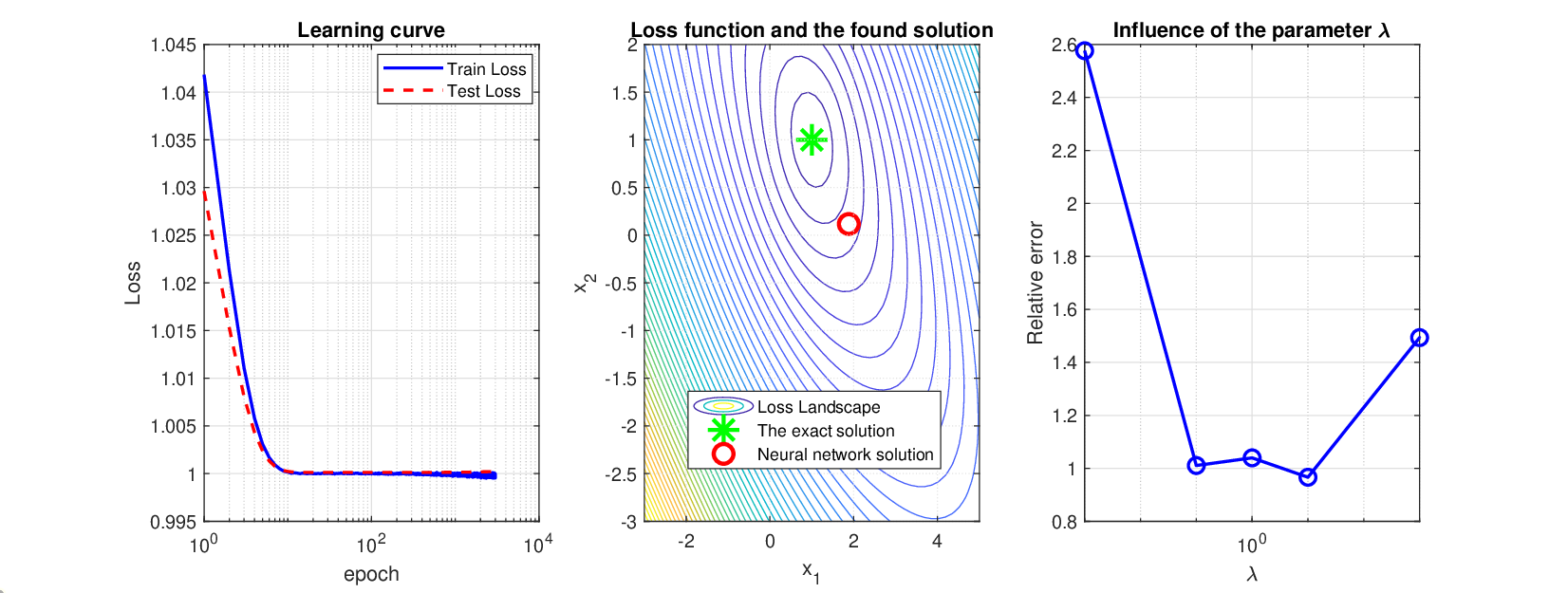}
\caption{Application of neural network regularization of weights for $h=0.0001$. Left: minimization of the loss function from \eqref{Loss1}; center: level lines of the loss function depending on unknowns $x_1, x_2$ at optimal weights and the optimal parameter $\lambda$; right: neural network accuracy as a function of the weight regularization parameter $\lambda$.}
\label{F5}
\end{figure}

The presented graphs demonstrate a significant advantage of the neural-network-based regularization algorithm over the conventional neural-network approach to solving this inverse problem. Although the NNBR algorithm solves an ill-posed problem with acceptable accuracy, albeit slightly lower than that of the standard TSVD regularization, a conventional neural network is incapable of yielding approximate solutions with acceptable accuracy. A similarly bad situation is observed in cases when the neural network was trained on a standard dataset using neural weights regularization at various levels of $\lambda$. An example of such regularization is shown in Fig. \ref{F5}.

\begin{figure}[t]
\centering
\includegraphics[width=1\textwidth]{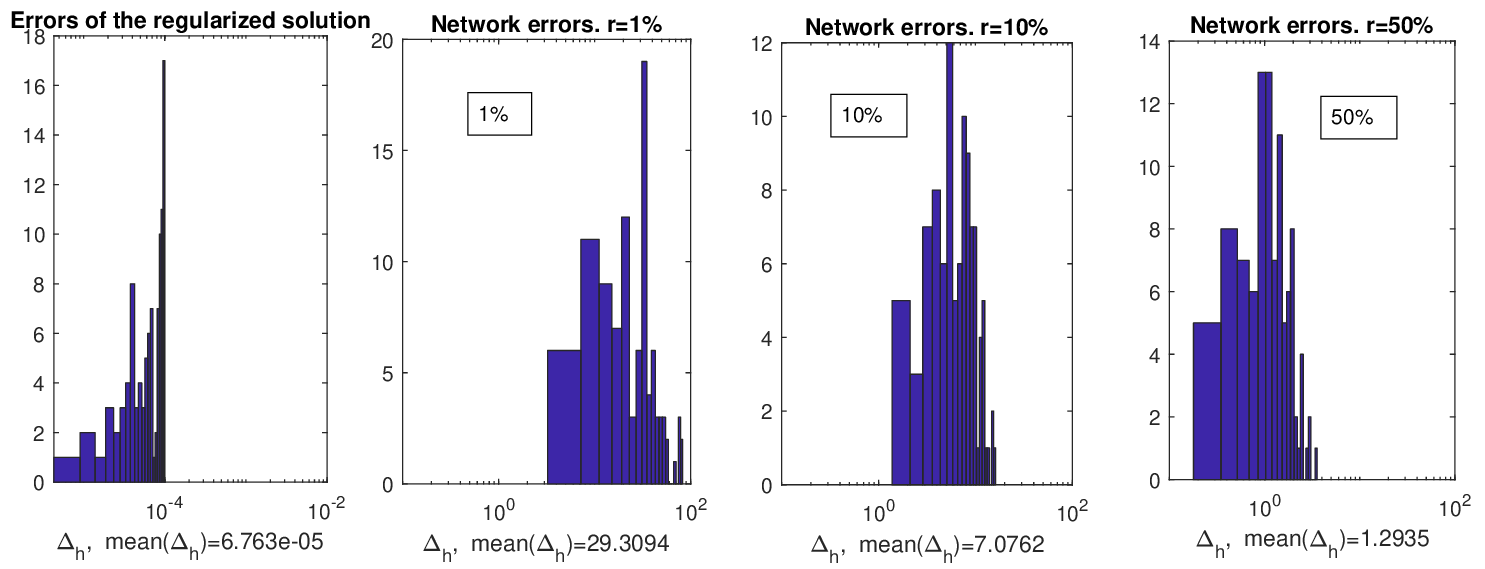}
\caption{Comparison of errors for the TSVD regularized solution and the standard neural networks solutions with dataset augmentation for $h=0.0001$. The datasets include data of the form $(A^+,A_h)$ in proportions of $r=1\%, 10\%, 50\%$. }
\label{F6}
\end{figure}

Another technique frequently employed to solve ill-posed problems using neural networks involves dataset augmentation (see, e.g., \cite{Aug}). This approach entails including in the dataset elements consisting of exact solutions and corresponding noisy data with a specified level of precision $\delta$: $(\bar x(y),{y_\delta })$. However, this approach does not guarantee the regularization of the solution. We demonstrate this in Fig. \ref{F6} for datasets that include the data $(A^+, A_h)$ for $h=0.0001$ in proportions of $r=1, 10, 50$ percent.
Similar results were obtained using neural network weight regularization for analogous datasets with augmentation.

\section{Conclusions}
1. The statements and calculation results presented above demonstrate a fundamental difference between classical regularization methods and well-known neural network algorithms for ill-posed problems. The latter are capable of solving only conditionally well-posed inverse problems, as they do not take into account Bakushinsky’s veto. Changes to the neural network architecture cannot alter this situation, as they merely modify the operator of the neural network solution to the problem but do not overcome Bakushinsky’s veto. Attempts to employ weight regularization in neural networks or to augment training sets by introducing perturbed data do not alter the situation from a general theoretical standpoint. Even for the simplest inverse problems, such as the one presented in Section \ref{Ill}, the errors associated with these ordinary neural network methods can be unacceptable.

2. Theoretically, an ill-posed inverse problem can be solved with the help of neural network methods only using a set of neural networks that explicitly depend on the known error levels $\delta$ of the approximate data. One way to achieve this is to use a regularized dataset constructed with a base regularization algorithm that explicitly depends on $\delta$. The resulting \emph{neural-network-based regularization algorithm} is capable of solving the ill-posed inverse problem. Test calculations demonstrate that the accuracy of the resulting approximate solutions is close in order of magnitude to that of the base regularization algorithm.

3. These conclusions hold true both for the application of abstract neural networks in infinite-dimensional spaces and for specific finite-dimensional neural networks. A typical example from Section \ref{Ill}, illustrating the use of neural-network-based approximation theory and ordinary neural network methods, confirms these conclusions.

\end{document}